\documentclass{article}
\usepackage{amssymb}
\usepackage{amsthm}
\usepackage{ifthen,amsfonts,graphicx,latexsym,algorithm,algorithmic,url}
\usepackage{amsmath,color}
\usepackage{verbatim}
\usepackage[applemac]{inputenc}
\usepackage[english]{babel}

\newtheorem{e-proposition}[theorem]{Proposition}

\theoremstyle{remark}

\newcommand{\pical}{\mathcal{P}}
\newcommand{\lcal}{\mathcal{L}}

\newcommand{\R}{\mathbb{R}}

\DeclareMathOperator*{\argmin}{argmin}

\def\og{\leavevmode\raise.3ex\hbox{$\scriptscriptstyle\langle\!\langle$~}}
\def\fg{\leavevmode\raise.3ex\hbox{~$\!\scriptscriptstyle\,\rangle\!\rangle$}}

\title{A Modest Proposal for MFG with Density Constraints\\{\small Mean Field Games and related topics, Rome, May 2011}}
\author{Filippo Santambrogio \thanks{Laboratoire de Mathématiques d'Orsay, Faculté de Sciences, Université Paris-Sud, 91405 Orsay cedex, FRANCE, {\tt filippo.santambrogio@math.u-psud.fr}}}

\begin{document}
\maketitle

\begin{abstract}
We consider a typical problem in Mean Field Games: the congestion case, where in the cost that agents optimize there is a penalization for passing through zones with high density of agents, in a deterministic framework. This equilibrium problem is known to be equivalent to the optimization of a global functional including an $L^p$ norm of the density. The question arises as to produce a similar model replacing the $L^p$ penalization with an $L^\infty$ constraint, but the simplest approaches do not give meaningful definitions. Taking into account recent works about crowd motion, where the density constraint $\rho\leq 1$ was treated in terms of projections of the velocity field onto the set of admissible velocity (with a constraint on the divergence) and a pressure field was introduced, we propose a definition and write a system of PDEs including the usual Hamilton-Jacobi equation coupled with the continuity equation. For this system, we analyze an example and propose some open problems.  \end{abstract}

\section{Introduction}

The goal of this paper, as well as that of the talk that the author gave at the workshop in Rome, is to present an idea to define a Mean Field Game (and to write down the corresponding coupled system of PDEs) when the moving density is subject to a congestion constraint, i.e. it cannot go beyond a fixed value standing for the maximal compression of the agents. We will see that this is not completely evident, but recalls what is done in some papers on ``hard'' congestion for crowd motion. In the same papers, the lack of strategical issues in the movement of the agents was perceived as a drawback, which motivates the definition of this new model. 

In the next section we will quickly review the existing models for non-strategic crowd motion and for deterministic Mean Field Games with density penalizations. Later on, in Section 3, we will explain how to glue these ideas so as to put constraints instead of penalizations and write down an MFG system. These systems usually couple a continuity equation and a Hamilton-Jacobi equation in the unknowns $\rho$ (density) and $\varphi$ (value function), but in this case a third unknown will appear, the pressure $p$, and the system will be complemented with some conditions on the pressure. Unfortunately, no result is available right now on this system, and Section 4 will be devoted to examples and open questions, in particular those arising from the comparison with the penalized models.

\section{About the existing models for collective motion under congestion effects}

We review in this section the models that are the closest to the one we look for. First we describe a model for crowd motion where the key point is the constraint $\rho\leq 1$, but no strategical issue or equilibrium condition is present; then we pass to a sketchy description of Mean Field Games with penalization of high densities. We will conclude with a naive, and wrong, idea to replace this penalization with a constraint. 

\subsection{Crowd motion with density constraints}

Many models for pedestrian motion are based on the idea that individual slow down when the density of the crowd is too high around them. This correspond to saying that their speed is (or is influenced by) a decreasing function of the density $\rho$. This is what is sometimes called ``soft congestion''. 

A more drastic viewpoint is the one presented by Maury and Venel in \cite{crowd1}, where the idea is that particles can move as they want as far as they are not too dense, but, if a density constraint is saturated, then their spontaneous velocity field $u$ could be impossible and would be abruptly modified into another field $v$, less concentrating and typically slower. The main assumption of such a model is that $v$ is the projection of $u$ onto the cone of admissible velocities, i.e. those that infinitesimally preserve the constraint. Maury and Venel were concerned with the discrete (microscopic) case, where individuals are represented by small disks. The density constraint is interpreted as a non-superposition constraint: particles cannot overlap, but as soon as they are not in touch their motion is unconstrained. When they touch, the set of admissible velocities restricts to those that increase the distance of every pair of particles in contact.

Later on, a  continuous (macroscopic) model has also been established in \cite{MauRouSan} (see also \cite{MauRouSanVen} for the latest developments and a micro-macro comparison).

  \begin{itemize}
  \item  The particles population is described by a probability density $\rho\in\pical(\Omega)$;
  \item the non-overlapping  constraint is replaced by the condition $\rho\in K=\{\rho\in\pical(\Omega)\,:\,\rho\leq 1\}$ (where $\rho$ denotes, by abuse of notation at the same time the probability and its density, since anyway we only consider absolutely continuous measures);
  \item for every time $t$, we consider $u_t:\Omega\to\R^d$ a vector field, possibly depending on $\rho$; 
  \item for every density $\rho$ we have a set of admissible velocities, characterized by the sign of the divergence on the saturated region $\{\rho=1\}$: $adm(\rho)=\big\{ v:\Omega\to\R^d\;:\;\nabla\cdot v\geq 0\;\mbox{ on }\{\rho=1\}\big\}$;
  \item we consider the projection $P$, which is either  the projection in $L^2(\lcal^d)$ or in $L^2(\rho)$ (this will turn out to be the same, since the only relevant zone is $\{\rho=1\}$); 
  \item we solve the equation 
  \begin{equation}\label{cont con proj}
 \partial_t\rho_t+\nabla\cdot\big( \rho_t \,\big(P_{adm(\rho_t)} [u_t]\big)\big)=0.
  \end{equation}
  \end{itemize}
  
Equation \eqref{cont con proj} is motivated by the fact that  the equation satisfied by the evolution of a density $\rho$ when each particle follows the velocity field $v$ is exactly the {\it continuity equation} $\partial_t\rho_t+\nabla\cdot\big( \rho_t v_t\big)=0$ (with $v\cdot n=0$ on $\partial\Omega$, so that the density does not exit $\Omega$). Here we only insert the fact that $v$ is the projection of $u$.
  
  The main difficulty is the fact that the vector field $v_t=P_{adm(\rho_t)} [u_t]$ is nor regular (since it is obtained as an $L^2$ projection, and may only be expected to be $L^2$ a priori), neither it depends regularly on $\rho$.
  
  Anyway, we need to make a bit more precise the definitions above. First, we stress that the continuity equation $\partial_t\rho+\nabla\cdot\big( \rho v\big)=0$ is to be intended in the weak sense:
  \begin{equation}\label{weakconf}
  \mbox{ for every }\psi\in C^1(\Omega),\quad t\mapsto\!\!\int_\Omega \psi d\rho_t\;\mbox{is Lipschitz and, for a.e. $t$ }\quad \frac{d}{dt}\int_\Omega \psi d\rho_t=\int_\Omega \nabla \psi \cdot v_t d\rho_t.
  \end{equation}
  Moreover, also the divergence of vector fields which are only supposed to be $L^2$ is to be defined in a weak sense, and to do that it is more convenient to give a better description of $adm(\rho)$ by duality: 
  $$adm(\rho)=\left\{ v\in L^2(\rho)\;:\;\int v\cdot\nabla p\leq 0\quad\forall p\in H^1(\Omega)\,:\,p\geq 0,\; p(1-\rho)=0\right\}.$$
 In this way we characterize $v=P_{adm(\rho)}[u]$ through 
   \begin{gather*}
  u=v+\nabla p,\quad v\in adm(\rho),\quad \int v\cdot\nabla p = 0,\\
  p\in press(\rho):=\{p\in H^1(\Omega),\, p\geq 0,\, p(1-\rho)=0\},
  \end{gather*}
  where $press(\rho)$ is the space of functions $p$ used as test functions in the dual definition of $adm(\rho)$. They play the role of pressures affecting the movement. The two cones $\nabla press(\rho)$ and $adm(\rho)$ are orthogonal cones and this allows for an orthogonal decomposition $u=v+\nabla p$. This also gives the alternative expression of Equation \eqref{cont con proj}, i.e.
  \begin{equation}\label{con press}
 \partial_t\rho+\nabla\cdot\big( \rho (u-\nabla p)\big)=0;\quad p\geq 0,\; p(1-\rho)=0.
  \end{equation}
We stress the fact that the orthogonality condition $\int v\cdot\nabla p=0$ is no more necessary and is a consequence of the continuity equation. This may be informally seen in the following way. Fix a time $t_0$ and notice that
$$\int p_{t_0}d\rho_{t_0}=\int p_{t_0}\geq \int p_{t_0}d\rho_{t}$$
which means that the function $t\mapsto \int p_{t_0}d\rho_{t}$ is maximal at $t=t_0$. By differentiating and using \eqref{weakconf}, we get 
$$0=\int \nabla p_{t_0}\cdot v_{t_0}d\rho_{t_0}=\int \nabla p_{t_0}\cdot v_{t_0}$$
(this proof is only formal because nothing guarantees that $t\mapsto \int p_{t_0}d\rho_{t}$ is differentiable at $t=t_0$, but this can be fixed and it is possible to obtain $\int \nabla p_{t}\cdot v_{t}=0$ for a.e. $t$).

As we said, this PDE is particularly difficult to analyze because of the lack of regularity of $v$, and the main tools to prove at least existence of a solution lie in the use of Optimal Transport and Wasserstein distances. We will not give here details on these topics, which are quite well known now and are not crucial for the subsequent parts of the paper. The interested reader may refer to \cite{AmbGigSav, Vil1, Vil2}. We will only describe here how Optimal Transport enters the game in the case where the spontaneous velocity $u$ is given by the gradient of a given function that agents want to minimize: $u=-\nabla D$.

Actually, it happens that in this case, Equation \eqref{cont con proj}, may be seen as the gradient flow of a functional $F$ with respect to the Wasserstein distance $W_2$ (which is the distance on the set $\pical(\Omega)$ induced by the minimal value of the optimal transport problem for the quadratic cost $|x-y|^2$). The functional that we need to consider is
$$F(\rho)=\begin{cases}\int D(x)d\rho & \mbox{ if }\rho\in K,\\
					+\infty & \mbox{ if }\rho\notin K,\end{cases}$$
where, again, $K=\{\rho\in\pical(\Omega)\,:\,\rho\leq 1\}.$ To explain what a gradient flow is, let us only say that it is the equation that rules an evolution where $\rho_0$ is given, and $\rho_t$ moves by selecting at every time the direction that lets the quantity $F(\rho_t)$ decrease as fast as possible. In particular this may be seen if one discretizes in time the problem, getting an  iterative method. We fix a time step $\tau>0$, take $\rho^\tau_0=\rho_0$ and at every step we find a new density through 
$$\rho^\tau_{k+1}\in\argmin_\rho F(\rho)+\frac{W_2^2(\rho,\rho^{\tau}_k)}{2\tau}=\argmin_{\rho\in K} \int D(x)d\rho+\frac{W_2^2(\rho,\rho^{\tau}_k)}{2\tau}.$$
Then, by doing it for smaller and smaller time steps $\tau$ and letting $\tau\to 0$, one recovers a curve of probabilities (continuous for the $W_2$ distance) and can prove that it is a solution of \eqref{cont con proj}. This follows the abstract ideas of De Giorgi and Ambrosio, see \cite{DeG,MovMin}, which are meant for general metric spaces (see also \cite{AmbGigSav}). When they are applied to the case of the Wasserstein distance they allow to obtain evolutionary PDEs, as it has been done first in \cite{JKO} for the Heat and Fokker-Plank equations.

An interesting point here is that the constraint $\rho\in K$ may be seen as a limit of $L^m$ penalizations as $m\to+\infty$. Indeed, one has
$$\lim_{m\to+\infty} \frac 1m \int \rho^m=\lim_{m\to+\infty} \frac 1m ||\rho||_{L^m}^m=\begin{cases}0&\mbox{ if }||\rho||_{L^\infty}\leq 1,\;\mbox{i.e.}\,\rho\in K\\
+\infty&\mbox{otherwise.}\end{cases}$$
This means that the functional $F$ may be approximated, as $m\to\infty$ with the functionals $F_m$ defined by $F_m(\rho)=\int D(x)\rho(x)dx+\frac 1m\int\rho(x)^mdx.$ The precise meaning of this approximation is to be read in terms of $\Gamma-$convergence (see \cite{introgammaconve}), and the applications of $\Gamma-$convergence to the corresponding gradient flows are presented in \cite{SerfatyGammaFlow}. This penalized counterpart will be the starting point for the next subsection. Here we only stress that the Gradient flow of the functional $F_m$ gives the well known porous media equations (see \cite{Ott2})
$$\partial_t\rho-\nabla\cdot\big( \rho \nabla D\big)-\frac{m-1}{m}\Delta\big(\rho^m\big)=0.$$

\subsection{Mean Field Games with congestion and density penalization}

Let us give a quick and informal presentation of the equations of Mean Field Games under congestion effects. As we said in the introduction, we will only consider the deterministic case (i.e. there will be no stochastic part in the evolution of the agents and we will use standard optimal control theory and not stochastic control). The reader may refer to \cite{LasLioCRAS, LasLio} and to the rest of this special issue; the presentation of the subject that is provided here does not want to be exhaustive but only to address the features of the problem that we need in the sequel.

Suppose that a population of agents may evolve in time, each agent following trajectories of the controlled equation
\begin{equation}\label{control}
y'(t)=f(t,y(y),\alpha(t)),\quad t\in [0,T]
\end{equation}
$\alpha:[0,T]\to\R^n$ being a control that every agent may choose. At each time $t$, the goal of each agent is to maximize the payoff
\begin{equation}\label{cost}
-\!\!\int_t^T \left(\frac{|\alpha(s)|^2}{2}+g(\rho(s,y(s)))\right)ds + \Phi(y(T)),
\end{equation}
where $g$ is a given increasing function. This means that $\alpha$ is the effort that every agent makes to move in the desired direction, and he pays for it (for simplicity, we take a quadratic cost), that its position depends on $\alpha$ through the equation \eqref{control}, and that he tries to optimize the final payoff $\Phi$ but he also pays for the densities of the regions he passes by. In particular agents would like to avoid crowded areas. At this first step the density $\rho(t,x)$ is supposed to be a given function. Yet, the MFG problem is an equilibrium problem: given the initial density $\rho_0$, find a time-dependent family of densities $\rho_t$ such that, when every agent selects its optimal trajectory so as to optimize \eqref{cost}, the density realized by these optimal trajectories at time $t$ are exactly $\rho_t$.

One can study the optimal control problem given by \eqref{control} and \eqref{cost} by means of its value function
$$\varphi(t,x)=\sup -\!\!\int_t^T \left(\frac{|\alpha(s)|^2}{2}+g(\rho(s,y(s)))\right)ds + \Phi(y(T))\;:\;y(t)=x,\;y'(s)=f(s,y(s),\alpha(s)).$$
It is well-known from optimal control theory that $\varphi$ satisfies the Hamilton-Jacobi-Bellmann equation (in the viscosity sense)
$$\partial_t \varphi(t,x) + H(t,x,\nabla\varphi(t,x))=0,\quad \varphi(T,x)=\Phi(x),$$
where the Hamiltonian $H$ is defined through
\begin{equation}\label{Hamiltonian}
H(t,x,\xi)=\sup_\alpha\; \xi\cdot f(t,x,\alpha)-\frac {|\alpha|^2}{2}-g(\rho_t(x)).
\end{equation}
Moreover, it is also well-known that in the control problem, for every $(t,x)$, the optimal choice of $\alpha(t)$ so as to optimize the criterion  \eqref{cost} starting from $x$ at time $t$ is the control $\alpha$ which maximizes in the definition of $H$ given in \eqref{Hamiltonian} for $\xi=\nabla \varphi (t,x)$, i.e. which maximizes $\nabla \varphi (t,x)\cdot f(t,x,\alpha)-\frac {|\alpha|^2}{2}$.

This gives a system of two coupled equations, since $\varphi$ solves an HJB equation where $\rho$ appears in the Hamiltonian, and $\rho$ evolves according to a continuity equation $\partial_t \rho + \nabla\cdot( \rho v)=0$, where the vector field $v(t,x)$ is given by $f(t,x,\alpha_{t,x})$ and the dependence of $\alpha_{t,x}$ on $(t,x)$ passes through $\nabla \varphi(t,x)$. To be more precise, we can give an explicit example, in the case $f(t,x,\alpha)=\alpha$ (this case corresponds to a classical calculus of variations problem, but we will still use the optimal control formalism since we will need to switch back to true control problems in a while). In such a case we can compute explicitly the Hamiltonian $H(t,x,\xi)$ thus getting $H(t,x,\xi)=\frac{|\xi|^2}{2}-g(\rho_t(x))$ and the optimal $\alpha$ in the definition of $H(t,x,\xi)$ turns out to be exactly $\xi$. Hence, the vector field to be put in the continuity equation is exactly $\nabla \varphi$. This gives the system
\begin{equation}\label{MFGgrho}
\begin{cases}\partial_t\varphi+\frac{|\nabla\varphi|^2}{2}-g(\rho)=0,\\
		\partial_t \rho +\nabla\cdot (\rho \nabla \varphi)=0,\\
		\varphi(T,x)=\Phi(x),\quad\rho(0,x)=\rho_0(x),\end{cases}
\end{equation}			
where the first equation is satisfied in the viscosity sense and the second in the distributional sense.

In this case, it is also known that a solution of this system (i.e. an equilibrium in the Mean Field Game) may be obtained by minimizing a suitable global functional (obviously, up to a change of sign, we could also express it as a maximization). Actually, one can solve
\begin{equation}\label{mini for mfg}
\min\quad \int_0^T\!\!\int_\Omega \left(\frac 12 |\alpha(t,x)|^2\rho(t,x)+G(\rho(t,x))\right)dxdt-\!\!\int_\Omega\Phi(x)\rho(T,x)dx
\end{equation}
among solutions $(\rho,\alpha)$ of the continuity equation $\partial_t \rho +\nabla\cdot (\rho \alpha)=0$ with initial datum $\rho(0,x)=\rho_0(x)$. When the function $G$ is chosen as the primitive $G$ of $g$ (i.e. $G'=g$, and in particular $G$ is convex), it happens that the minimizers of this global functional are equilibria in the sense explained above (the trajectories that the agents chose according to these densities are such that they exactly realize these densities). These functionals recall the functional proposed by Benamou and Brenier \cite{BenBre} to give a dynamical formulation of optimal transport, or those studied by Buttazzo, Jimenez and Oudet for density congested dynamics in \cite{ButJimOud}.

It is clear that we can choose $g(\rho)= \rho^{m-1}$ (for $m>1$), which gives $G(\rho)=\frac 1m \rho^m$, and the minimization problem in \eqref{mini for mfg} tends, as $m\to+\infty$ to  
\begin{equation}\label{mini for mfg Linfty}
\min\quad \int_0^T\!\!\int_\Omega\frac 12 |\alpha(t,x)|^2\rho(t,x) dxdt-\!\!\int_\Omega\Phi(x)\rho(T,x)dx,
\end{equation}
the minimization being, again, performed among solutions $(\rho,\alpha)$ of the continuity equation $\partial_t \rho +\nabla\cdot (\rho \alpha)=0$ with initial datum $\rho(0,x)=\rho_0(x)$ and satisfying the constraint $\rho(t,x)\leq 1$ a.e. In particular, the energy is not finite unless $\rho_0$ satisfies the constraint $\rho_0\leq 1$ as well. By the way, this problem is equivalent, at least in a convex domain $\Omega$, to finding $\rho_T\leq 1$ so as to minimize $T W_2^2(\rho_0,\rho_T)-\int \Phi d\rho_T$ and then choosing $\rho_t$ as a constant speed geodesic connecting $\rho_0$ to $\rho_T$. This comes from the fact that the optimal curve connecting $\rho_0$ to $\rho_T$ and minimizing $\int\!\!\int \frac 12 |\alpha|^2\rho dxdt$ is the geodesic between these two measures in the Wasserstein space, and the set $K$ of measures with density bounded by $1$ is known to be geodesically convex (i.e., if two measures belong to $K$, the geodesic between them stays in $K$ as well). This concept of convexity and the proof that $K$ is geodesically convex come from \cite{MC}.

The question arises naturally of how to define a Mean Field Game equilibrium problem under the density constraint $\rho\leq 1$ instead of the density penalization $+\rho^{m-1}$.

It happens that this question is not trivial at all. Let us try the simplest definition. We say that a curve of densities $\rho_t$ is an equilibrium if 
\begin{itemize}
\item it satisfies the constraint $||\rho_t||_{L^\infty}\leq 1$ for all times $t$,
\item when the agents choose their optimal trajectory, solving  
$$
\max\quad -\!\!\int_t^T \frac{|\alpha(s)|^2}{2}ds + \Phi(y(T)),\;:\;y'(s)=\alpha(s),\;y(t)=x
$$
and respecting the constraint $\rho_s(y(s))\leq 1$, the density they realize at time $t$ is again $\rho_t$.
\end{itemize}
Yet, this definition does not make any sense (and would lead in general to non-existence of such an equilibrium) because of this fact:  this game is a non-atomic game, where the influence of every individual is negligible; thus, every agent facing densities $\rho_t\leq 1$ can choose with no restrictions at all its trajectory, since he knows that anyway simply adding its own presence somewhere will not really change the density and violate the constraint. This means that the trajectories that are selected by every agent do not take into account the constraint $\rho\leq 1$ and hence there is no reason for the densities realized by these choices to satisfy $\rho\leq 1$.

\section{The MFG system with density constraint}

The main idea of this section, and of this ``modest proposal'', is the following: use the pressure! Express the effects of the density constraints through the existence of a non-zero pressure field.

This requires changing the definition of the equilibrium.
\begin{itemize}
\item We describe the situation through a pair $(\rho,\bar\alpha)$ where $\rho_t$ stands for the density of the agents at time $t$ and $\bar\alpha(t,x)$ for the effort that the agent located at $x$ at time $t$ makes to control his movement. We require that $(\rho,\bar\alpha)$ satisfies
$$\partial_t\rho_t+\nabla\cdot\big( \rho_t \,\big(P_{adm(\rho_t)} [\bar\alpha_t]\big)\big)=0,$$
which means that $\rho$ is advected by the projection of this effort vector field.
\item Considering the projection of $\bar\alpha_t$ onto $adm(\rho_t)$ a pressure $p_t$ appears such that
$$P_{adm(\rho_t)}[\bar \alpha_t]=\bar\alpha_t-\nabla p_t,\quad p_t\in press(\rho_t),\;\int P_{adm(\rho_t)} [\alpha_t]\cdot \nabla p_t=0$$
(but, again, the orthogonality condition is included in the continuity equation).
\item Every agent wants to optimize his own control problem where the state equation is influenced by the pressure:
$$
\max\quad -\!\!\int_t^T \frac{|\alpha(s)|^2}{2}ds + \Phi(y(T)),\;:\;y'(s)=\alpha(s)-\nabla p_s(y(s)),\;y(t)=x.
$$
\item The configuration $(\rho,\bar\alpha)$ is said to be an equilibrium if the original effort field $\bar\alpha$ coincides with the optimal effort in this control problem and if the original densities $\rho_t$ coincide with the densities realized at time $t$ according to these optimal trajectories.
\end{itemize}

In order to write down the equations that correspond to such an equilibrium we first need to compute the Hamiltonian of this control problem:
$$H(t,x,\xi)=\sup_\alpha\quad \xi\cdot(\alpha-\nabla p(t,x))-\frac{|\alpha|^2}{2}\;=\;\frac{|\xi|^2}{2}-\xi\cdot\nabla p(t,x),$$
the optimal $\alpha$ in this maximization being exactly $\xi$. This means that the optimal effort at $(t,x)$ will be $\nabla\varphi(t,x)$ and that the vector field appearing in the continuity equation will be  $\nabla\varphi-\nabla p$.
We finally get to the MFG system
\begin{equation}\label{MFGdensity1}
\begin{cases}\partial_t\varphi+\frac{|\nabla\varphi|^2}{2}-\nabla\varphi\cdot\nabla p=0,\\
		\partial_t \rho +\nabla\cdot (\rho (\nabla \varphi-\nabla p))=0,\\
		p\geq 0,\, p(1-\rho)=0,\\
		\varphi(T,x)=\Phi(x),\quad\rho(0,x)=\rho_0(x).\end{cases}
\end{equation}

It is important to notice that in this case the equilibrium may only be defined in terms of a pair $(\rho,\bar\alpha)$, and that the density only is not sufficient to describe the configuration. The reason lies in the fact that the pressure at time $t$ does not depend on the density $\rho_t$ only, but on the vector field that we project on $adm(\rho_t)$. This was not the case in the gradient flow framework of crowd motion in \cite{MauRouSan}, since the vector field to be projected was itself a function of $\rho$. This will be clearer in the next example.

\section{An example and some open questions}

\subsection{An example where nothing moves}
 
 Even if the example presented in this section is very simple, its resolution will be quite sketchy and informal.
 Consider the following situation: given the function $\Phi$, suppose that there exists $\ell\in\R$ such that $|\{\Phi>\ell\}|=1$. We will also suppose $\Phi$ to be regular enough (in particular, to be continuous), so that $\Phi=\ell$ on $\partial A$. Let us set $A=\{\Phi>\ell\}$ and take $\rho_0=I_A$, i.e. $\rho_0(x)=1$ for $x\in A$ and $\rho_0(x)=0$ for $x\notin A$.
 
 Let us try to guess the solution of  \eqref{MFGdensity1}. If instead of a density constraint we had a power penalization, the solution would be obtained by minimization of a global functional $\int\!\!\int\rho\alpha^2+\frac 1m\int\!\!\int \rho^m-\!\!\int\Phi d\rho_T$. Let us suppose for a while that a similar principle holds for the $L^\infty$ constraint: we could imagine that the solution is obtained by minimizing  $\int\!\!\int\rho\alpha^2-\!\!\int\Phi d\rho_T$ under the constraint $\rho\leq 1$. This problem has an easy solution: not moving at all! Actually, $\rho_t=\rho_0$ is a solution since it minimizes the first term (no movement, no kinetic energy) and maximizes $\int\Phi d\rho$ among  probability measures with $\rho\leq 1$, since none of them can do better than $\rho_0$, which is already concentrated on the best points for the function $\Phi$. This suggests that the solution of \eqref{MFGdensity1} should be obtained with $\rho_t=\rho_0$.
 
 Can we conclude that the pair $(\rho,\bar\alpha)$ given by $\rho_t=\rho_0$ and $\bar\alpha_t=0$ is an equilibrium? the answer is for sure not. Why? just because if the effort vector field is identically zero, then its projection will also be zero, and the pressure as well. This means that we would have $\nabla p_t=0$, i.e. the agents could move as if there were no constraints. Hence, every agent would maximize his own payoff $-\!\!\int_t^T \frac{|\alpha(s)|^2}{2}ds + \Phi(y(T)),\;:\;y(t)=x,\;y'(s)=\alpha(s),$ and in general the solution is not $\alpha=0$.
 
 Yet, even if we did not find an equilibrium so far, we can improve our guess. We can claim that the equilibrium configuration is such that everybody wants to move, but nobody manages to do it. This means $\rho_t=\rho_0$ but $\bar\alpha\neq 0$. It is indeed possible that the optimal effort $\bar\alpha$ is exactly compensated by the pressure effects. Let us look for a solution of \eqref{MFGdensity1} of this form. 
 
 We need to impose $\nabla\varphi=\nabla p$ on $[0,T]\times A$. This is sufficient to satisfy the continuity equation with a static measure $\rho_t=\rho0$. If we consider that $p$ must vanish where $\rho<1$, this imposes $p=0$ on $[0,T]\times A^c$. Hence we get the (strange) equation satisfied (at least formally) by the value function $\varphi$:
 \begin{equation}\label{phipm}
 \partial_t\varphi+\frac{|\nabla\varphi|^2}{2}=0 \mbox{ on }[0,T]\times A^c\quad\mbox{ and }\quad\partial_t\varphi-\frac{|\nabla\varphi|^2}{2}=0 \mbox{ on }[0,T]\times A.
 \end{equation}
 We need to solve this equation complemented by its final condition $\varphi(T,x)=\Phi(x)$. Once we find $\varphi$, using $\nabla\varphi=\nabla p$ on $[0,T]\times A$, we need to deduce a pressure $p$ and to prove that it is actually a pressure, i.e. $p\geq 0$ and $p(1-\rho)=0$.
 
 An informal way to solve \eqref{phipm} is the following. For simplicity, let us suppose $\ell=0$ (which is always possible up to adding a constant to $\Phi$). Consider the function $\tilde\Phi=-|\Phi|$ and define
 $$\tilde\varphi(t,x)=\sup -\!\!\int_t^T\frac{|\alpha(s)|^2}{2}ds + \tilde\Phi(y(T))\;:\;y(t)=x,\;y'(s)=\alpha(s).$$
Let us suppose that this value function is regular: it must satisfy 
$$\partial_t\tilde\varphi+\frac{|\nabla\tilde\varphi|^2}{2}=0,\quad\tilde\varphi(T,x)=\tilde\Phi(x),\quad\tilde\varphi\leq 0,\quad\tilde\varphi(t,x)=0\mbox{ for } x\in \partial A.$$
The last conditions are easy to check: we see from the definition and from $\tilde\Phi\leq 0$ that we have $\tilde\varphi\leq 0$; yet, the constant curve $y(s)=x$ gives $\tilde\varphi(t,x)=0$, provided $\tilde\Phi(x)=0$, which is the case for $x\in\partial A$. After that one can define 
$$\varphi=(1-2I_A)\tilde\varphi,$$
which means $\varphi=\tilde\varphi$ outside $A$ and $\varphi=-\tilde\varphi$ on $A$. It is possible to check that (at least in a formal way) $\varphi$ solves \eqref{phipm} and $p$ is an admissible pressure.

Obviously in the analysis of this example we just provided a solution, with no guarantee that it is the unique one.

\subsection{Some open questions}

As we said in the introduction, in this paper we only present System \eqref{MFGdensity1}, its motivations and connections with other problems, but we are not able to prove any result on it. In particular existence, uniqueness and well-posedness in general are open questions.

The usual conditions guaranteeing uniqueness in MFG are based on sign conditions such as
$$\frac{d}{dt}\int(\varphi_1-\varphi_2)d(\rho_1-\rho_2)\geq 0,$$
to be checked on every pair of solutions $(\varphi_1,\rho_1),\,(\varphi_2,\rho_2)$. Yet, they do not seem easy to verify here, nor to imply in a simple way uniqueness results, even if some estimates are possible.

Hence, we ignore for the moment these uniqueness questions, and we only concentrate on two questions concerning the connection with the $\rho^m$ penalization and a possible variational principle.

\paragraph{First question: can we obtain \eqref{MFGdensity1} as a limit of $L^m$ penalizations with $m\to\infty$?} Let us be more precise. Consider the System \eqref{MFGgrho}, choosing $g(\rho)=\rho^{m-1}$ (which corresponds to a penalization $\frac 1m \int\rho^m$ in the global functional that is optimized by the equilibrium). Obviously, in this system there is no pressure, but the question is: is it possible to introduce fictitiously a function playing the role of the pressure, thus obtaining a 3-tuple $(\rho_m,\hat\varphi_m,\hat p_m)$ solving System \eqref{MFGgrho}, and then wonder if they converge, as $m\to\infty$, to a solution $(\rho,\varphi,p)$ of \eqref{MFGdensity1}?

To make an example of what we mean, we can decide to write the value function $\varphi$ in \eqref{MFGgrho} as $\varphi=\hat\varphi-\hat p$. We will choose later how to define $p$. If we rewrite the first equation of the system we have
$$\partial_t\hat\varphi-\partial_t\hat p +\frac 12|\nabla\hat\varphi|^2-\nabla\hat\varphi\cdot\nabla\hat p +\frac 12|\nabla\hat p|^2-\rho^{m-1}=0.$$
This means that, if we choose $\hat p$ satisfying
$$\begin{cases}-\partial_t\hat p+\frac 12|\nabla\hat p|^2-\rho^{m-1}=0,\\
			\hat p(T,x)=0,\end{cases}$$
the two equations of the system become exactly the same as those of \eqref{MFGdensity1} and $\hat\varphi(T,x)=\Phi(x)$. It is natural to wonder whether the function $\hat p$ that we just defined has some chances to converge towards the pressure $p$ of the limit system, and in particular if it will satisfy, at the limit, the condition $\hat p=0$ on $\{\rho<1\}$. This is not evident, and the only reasonable tool seems to interpret the condition $-\partial_t\hat p+\frac 12|\nabla\hat p|^2-\rho^{m-1}=0$ as a Hamilton-Jacobi equation for a control problem, thus writing
$$\hat p(t,x)=\inf\int_t^T\left(\frac{|\alpha(s)|^2}{2}+\rho^{m-1}(y(s))\right)ds \;:\;y(t)=x,\;y'(s)=\alpha(s).$$

\paragraph{Second question: is the equilibrium optimal in some sense?} We know that for the $L^m$ penalization the equilibrium is also obtained if one optimizes the global cost \eqref{mini for mfg}. One could think that letting $m\to+\infty$, the constrained MFG equilibrium should solve
\begin{equation}\label{minrholeq}
\min \quad \int_0^T\!\!\int_\Omega \left(\frac 12 |v(t,x)|^2\rho(t,x)\right)dxdt-\int_\Omega\Phi(x)\rho(T,x)dx
\end{equation}
among solutions $(\rho,v)$ of the continuity equation $\partial_t \rho +\nabla\cdot (\rho v)=0$ with initial datum $\rho(0,x)=\rho_0(x)$, satisfying the constraint $\rho\leq 1$. First, it is important to notice that here we consider the true velocity field $v$ and not the effort field $\alpha$, and we know that the equilibrium in this constrained setting cannot be expressed in terms of $(\rho_t,v_t)$ only. This means that, if this variational principle is true, then one could use it to identify the velocity field $v$, but would still need to decompose it into $\nabla\varphi-\nabla p$. Yet, it would allow to find the curve of densities $\rho_t$. 

It is an open question whether this variational principle holds or not. It holds in the example of the previous section. To check it in full generality, one possibility would be to use the equations of \eqref{MFGdensity1} to check whether the necessary conditions for optimality are satisfied.

Actually, one can write down by duality some necessary conditions for solving \eqref{minrholeq} under the constraint $\rho\leq 1$. In order to apply duality, it is well known, from Benamou-Brenier on (\cite{BenBre}), that it is convenient to use the variables $(\rho,q)$ with $q=\rho v$ instead of $(\rho,v)$. This leads to the optimization problem
$$\min_{(\rho,q)\,:\,0\leq\rho\leq 1}\int\!\!\int\frac{|q|^2}{2\rho}-\!\!\int\Phi d\rho_T+\sup_\chi -\!\!\int\!\!\int\partial_t \chi d\rho-\!\!\int\!\!\int\nabla \chi dq+\int\chi_Td\rho_T-\!\!\int\chi_0d\rho_0,$$
where the supremum over $\chi$ stands for the constraint given by the continuity equation. It may be dualized by interchanging $\min$ and $\sup$ and, optimizing first over $q$, one finds that the optimal $q$ is given by $q=\rho\nabla\chi$ and gets to
$$\sup_\chi\min_{0\leq\rho\leq 1} \int(\chi_T-\Phi)d\rho_T -\!\!\int \chi_0d\rho_0-\!\!\int\!\!\int\left(\partial_t \chi+\frac 12 |\nabla\chi|^2\right) d\rho.$$
This allows (we skip all justifications coming from convex analysis) to write the following optimality conditions : if $(\rho,q)$ is optimal then there exists a function $\chi$ such that $q=\rho\nabla\chi$ (and hence $v=\nabla\chi$) and
\begin{itemize}
\item $\partial_t \chi+\frac 12 |\nabla\chi|^2\leq 0$ on $\{\rho=0\}$ (but this is not restrictive, since $\chi$ may be modified on $\rho$-negligible sets and still satisfy $q=\rho\nabla\chi$ ),
\item  $\partial_t \chi+\frac 12 |\nabla\chi|^2\geq 0$ on $\{\rho=1\}$,
\item  $\partial_t \chi+\frac 12 |\nabla\chi|^2= 0$ on $\{0<\rho<1\}$,
 \item $\chi_T-\Phi\geq 0$ on $\{\rho_T=0\}$ (same consideration: it is not restrictive),
 \item $\chi_T-\Phi\leq 0$ on $\{\rho_T=1\}$, 
 \item $\chi_T-\Phi= 0$ on $\{0<\rho_T<1\}$.
 \end{itemize}
 
 The question of the optimality of the solution of \eqref{MFGdensity1} becomes now a question on the conditions satisfied by its velocity field. Since we already know $v=\nabla\varphi-\nabla p$, it is sufficient to define $\chi=\varphi-p$ and investigate whether the above conditions are satisfied or not.
 In particular we have
 $$\partial_t \chi+\frac 12 |\nabla\chi|^2=\partial_t \varphi+\frac 12 |\nabla\varphi|^2-\nabla\varphi\cdot\nabla p-\partial_t p+\frac 12 |\nabla p|^2=-\partial_t p+\frac 12 |\nabla p|^2,$$
 and the question is whether this quantity is positive on $\{\rho=1\}$ and zero on $\{0<\rho<1\}$. This question seems linked to the previous one, where the limit $m\to+\infty$ let appear exactly the condition $=-\partial_t \hat p+\frac 12 |\nabla\hat p|^2=\rho^{m-1}$.
 
 The conditions to be verified at time $t=T$ are, instead, easy to check, since $(\chi_T-\Phi)(x)=\varphi(T,x)-p(T,x)-\Phi(x)=-p(T,x)$ and we know $p\geq 0$ and $p=0$ on $\{0<\rho<1\}$.

 \paragraph{Acknowledgments}
 
 The author wants to thank the organizers of the Workshop ``Mean Field Games and related topics'' for the invitation to give a talk on a subject he is not a specialist of. For the same reason he is also indebted to several interesting discussions with several colleagues, and in particular G. Carlier and P. Cardaliaguet. The support of the ANR project EVaMEF ANR-09-JCJC-0096-01 is also acknowledged. Finally, the author thanks J. Swift for the first words of the title.

\end{document}